\documentclass[12pt]{article}
\usepackage[english]{babel}
\usepackage{amssymb}
\usepackage{inputenc}
\begin{document}

\author{Sh. A. Ayupov $^{1},$ A. A. Zaitov $^2$}

\title{\bf On the weight and the density of the space of order-preserving
functionals}

\maketitle

\begin{abstract}
In the present paper it is proved that the functors $O_\tau$ of
$\tau$-smooth order preserving functionals and $O_R$ of Radon order
preserving functionals preserve the weight of infinite Tychonoff
spaces. Moreover, it is established that the density and the weak
density of infinite Tychonoff spaces do not increase under these
functors.
\end{abstract}

\medskip

$^1$ Institute of Mathematics and information technologies,
Uzbekistan Academy of Science, F. Hodjaev str. 29, 100125, Tashkent
(Uzbekistan), e-mail: \emph{sh\_ayupov@mail.ru,
e\_ayupov@hotmail.com, mathinst@uzsci.net}

 $^{2}$ Institute of Mathematics and information technologies, Uzbekistan Academy of
Science, F. Hodjaev str. 29, 100125, Tashkent (Uzbekistan), e-mail:
\emph{adilbek$_{-}$zaitov@mail.ru}



\newpage
\large

\section*{\center 0. Introduction}

Let $X$ be a compact($\equiv$compact Hausdorff topological space)
and let $C(X)$ be the Banach algebra of all continuous real-valued
functions with the usual algebraic operations and with the
$\sup$-norm. For functions $\varphi,\ \psi\in C(X)$ we shall write
$\varphi\leq\psi$ if $\varphi(x)\leq\psi(x)$ for all $x\in X.$ If
$c\in \mathbb{R}$ then by $c_X$ we denote the constant function
identically equal to $c.$ Recall that a functional $\mu :
C(X)\rightarrow \mathbb{R}$ is said [2] to be:

1) \textit{order-preserving} if for any pair $\varphi, \psi\in C(X)$
of functions the inequality $\varphi\leq\psi$ implies
$\mu(\varphi)\leq\mu(\psi)$;

2) \textit{weakly additive} if
$\mu(\varphi+c_X)=\mu(\varphi)+c\mu(1_X)$ for all $\varphi\in C(X)$
and $c\in R$;

3) \textit{normed} if $\mu(1_X)=1$.

For a compact $X$ denote by $O(X)$ the set of all order-preserving
weakly additive and normed functionals $\mu : C(X)\rightarrow
\mathbb{R}$. By $W(X)$ we denote the set of all functionals
satisfying only the conditions 1) and 2) of the above definition.
Note that according to proposition 1 [7] each  order-preserving
weakly additive functional is continuous. Further order-preserving
weakly additive functionals are called order-preserving functionals
[2].

Let $X$ be a Tychonoff space and let $C_b(X)$ be the algebra of all
bounded continuous real-valued functions with the pointwise
algebraic operations. For a function $\varphi\in C_b(X)$ put
$\|\varphi\|=\sup\{|\varphi(x)|:\ x\in X\}.$ $C_b(X)$ with this norm
is a Banach algebra. For a net $\{\varphi_\alpha\}\subset C_b(X)$
$\varphi_\alpha \downarrow 0_X$ means that for every point $x\in X$
one has $\varphi_\alpha(x) \geq \varphi_\beta(x)$ at $\beta \succ
\alpha$ and $\lim\limits_\alpha \varphi_\alpha(x)=0_X.$ In this case
we say that $\{\varphi_\alpha\}$ is a monotone decreasing net
pointwise convergent to zero.

For a Tychonoff space $X$ by $\beta X$ denote its Stone-\v{C}ech
compact extension. Given any function $\varphi\in C_b(X)$ consider
its continuous extension $\widetilde{\varphi}\in C(\beta X).$ This
gives an isomorphism between the spaces $C_b(X)$ and $C(\beta X)$
moreover, $\|\widetilde{\varphi}\|=\|\varphi\|,$ i. e. this
isomorphism is an isometry, and topological properties of the above
spaces coincide. Therefore one may consider any function from
$C_b(X)$ as an element of $C(\beta X).$ Hence definitions 0.1 and
0.2 from [1] may be given in the following form.

\textbf{Definition 1.} An order-preserving functional $\mu\in
W(\beta X)$ is said to be \textit{$\tau$-smooth} if $\mu(\varphi
_\alpha)$ $\rightarrow 0$ for each monotone  net $\{\varphi_\alpha\}
\subset C(\beta X)$ decreasing to zero on $X.$

 \textbf{Definition 2.} An order-preserving functional $\mu\in O(\beta X)$ is said to be
\textit{Radon} order-preserving functional if $\mu(\varphi_\alpha)$
$\rightarrow 0$ for each bounded net $\{\varphi_\alpha \} \subset
C(\beta X)$ which uniformly converges to zero on compact subsets of
$X$.

For a Tychonoff space $X$ by $W_\tau(X)$ and $W_R(X)$ denote the
sets of all $\tau$-smooth and Radon order-preserving functionals
from $W(\beta X),$ respectively. The sets $W_\tau(X)$ and $W_R(X)$
are equipped with the pointwise convergence topology. The base of
neighborhoods of a functional $\mu \in W_\tau (X)$ (respectively, of
$\mu\in W_R(X)$) in the pointwise convergence topology consists of
the sets
$$
\langle \mu; \varphi_1,...,\varphi_k; \varepsilon \rangle = \{\nu
\in W(\beta X): |\nu(\varphi_i) - \mu(\varphi_i)| < \varepsilon \}
\cap W_\tau (X)
$$
$$
(\mbox{respectively, }\langle \mu; \varphi_1,...,\varphi_k;
\varepsilon \rangle = \{\nu \in W(\beta X): |\nu(\varphi_i) -
\mu(\varphi_i)| < \varepsilon \}\cap W_R(X))
$$
where $\varphi_i \in C(\beta X),$ $i=1,...,k$ and $\varepsilon > 0$.

Put
$$
O_\tau(X)=\{\mu\in W_\tau(X): \mu(1_X)=1\},
$$
$$
O_R(X)=\{\mu\in W_R(X): \mu(1_X)=1\}.
$$

The operations $O_\tau$ and $O_R$ are functors [1, Theorem 0.3] in
the category $Tych$ of Tychonoff spaces and their continuous maps.

Let $A$ be a closed subset of the compact $X.$ An order-preserving
functional $\mu\in O(X)$ is said to be \textit{supported} on $A$ if
$\mu\in O(A)$ [2]. The set
$$
\mbox {supp}\mu=\cap\{A:\ \mu\in O(A) \mbox{ and } A \mbox{ is
closed in } X\}
$$
is called the \textit{support} of the order-preserving functional
$\mu.$

For a Tychonoff space $X$ put
$$
O_\beta(X)=\{\mu\in O(\beta X): \mbox {supp} \mu\subset X\}.
$$
The operation $O_\beta$ translating a Tychonoff space $X$ to
$O_\beta(X),$ is a functor [3] in the category $Tych$. Obviously the
inclusions
$$
O_\beta (X) \subset O_R (X) \subset O_\tau (X) \subset O(\beta X)
\eqno (1)
$$
are valid for any Tychonoff space $X,$ and the equalities
$$
O_\beta (X) = O_R (X) = O_\tau (X) = O(\beta X)
$$
are true for arbitrary compact $X.$

Let $X$ and $Y$ be compacts and let $f:X\rightarrow Y$ be a
continuous map. Then the map $O(f):O(X)\rightarrow O(Y)$ defined by
the formula
$$
O(f)(\mu)(\varphi)=\mu(\varphi \circ f)
$$
is continuous, where $\varphi\in C(Y)$ and $\mu\in O(X).$

Now let $X$ and $Y$ be Tychonoff spaces and let $f:X\rightarrow Y$
be a continuous map. Put
$$
O_\tau(f)=O(\beta f)|O_\tau(X),
$$
$$
O_R(f)=O(\beta f)|O_R(X)
$$
and
$$
O_\beta(f)=O(\beta f)|O_\beta(X)
$$
where $\beta f: \beta X\rightarrow\beta Y$ is the Stone-\v{C}ech
extension of $f.$

Note that the above maps $O_\tau(f):O_\tau(X)\rightarrow O_\tau(Y),$
$O_R(f):O_R(X)\rightarrow O_R(Y)$ and
$O_\beta(f):O_\beta(X)\rightarrow O_\beta(Y)$ are defined correctly
and they also are continuous.

In the paper [3] it was shown that the functor $O_\beta$ preserves
the weight of infinite Tychonoff spaces. There exists an example [4,
example 3], which shows that under the functor $O(\beta\cdot)$ the
weight of a Tychonoff space may strictly increase, more precisely,
the example shows that the functor translates a Tychonoff space $X$
with countable weight to a compact $O(\beta X)$ with continuum
weight. The question whether the functors $O_\tau$ and $O_R$
preserve the weight was open. In this paper we obtain a positive
answer to this question. Moreover we prove that under the functors
$O_\tau$ and $O_R$ the density and the weak density of infinite
Tychonoff spaces do not increase. It is established that the weak
densities of the following spaces coincide:

$O_\omega(X)$ of order-preserving functionals with finite supports,

$O_\beta(X)$ of order-preserving functionals with compact supports,

$O_R(X)$ of Radon order-preserving functionals,

$O_\tau (X)$ of $\tau$-smooth order-preserving functionals and

$O(\beta X)$ of all order-preserving functionals.

Note that the space $W_\tau(X)$ equipped with the pointwise
convergence topology  may be considered as a subspace of the
topological product $\Pi=\prod\{\mathbb{R}_\varphi: \varphi\in
C(\beta X)\}$ of real lines $\mathbb{R}_\varphi=\mathbb{R}$. Since
$\Pi$ is a Tychonoff space the spaces $W_\tau(X)$ and $O_\tau(X)$
with the pointwise convergence topology  are also Tychonoff spaces.

\section*{\center 1. Main results}

Let $X$ be a topological space. Recall that a \textit{weight} of $X$
is the cardinal number $w(X)$ defined by the formula
$$
w(X)=\min\{|\mathfrak{B}|:\mathfrak{B} \mbox{ is a base of the
topology on } X\}.
$$
In this section we shall prove that the functors $O_\tau$ of
$\tau$-smooth order-preserving functionals and $O_R$ of Radon
order-preserving functionals preserve the weight of infinite
Tychonoff spaces. To do this, we need some constructions.

Let $Y$ be a subspace of a Tychonoff space $X.$ Put
$$
C=\{\psi|Y:\psi\in C_b(X)\}.
$$

The following notion is well-known. A subspace $Y\subset X$ is
called \textit{$C$-embedded} in $X$ if for each function $\varphi\in
C_b(Y)$ there exists a function $\widetilde{\varphi}\in C_b(X)$
such, that $\widetilde{\varphi}|Y=\varphi.$ If $Y$ is a $C$-embedded
subspace of the given space $X$ then clearly $C\equiv C_b(Y).$

For a $C$-embedded subspace $Y$ of a Tychonoff space $X,$ a
functional $\mu\in W(\beta X)$ and a function $\varphi\in C_b(Y)$
put
$$
r^X_Y(\mu)(\varphi)=
$$
$$
=\inf\{\mu(\psi):\psi\in C_b(X),\ \psi\geq (inf \{\varphi(y): y\in
Y\})_X,\ \psi|Y=\varphi\}. \eqno (2)
$$

\textbf{Lemma 1.} \textit{For each $\mu\in W(\beta X)$ we have
$r^X_Y(\mu)\in W(\beta Y).$ In other words $r^X_Y(\mu)$ is an
order-preserving weakly additive functional on $C_b(Y).$}

\textsl{Proof.} We have the following equality
$$
r^X_Y(\mu)(1_Y)=\mu(1_X), \eqno (3)
$$
which directly follows from (2). Let $\varphi\in C_b(Y).$ Then
$\varphi+c_Y\in C$ for all $c\in \mathbb{R}.$ We have
$$
r^X_Y(\mu)(\varphi+c_Y)=
$$
$$
=inf\{\mu(\psi):\psi\in C_b(X), \psi\geq(inf\{\varphi(y)+c: y\in
Y\})_X, \psi|Y = \varphi+c_Y\} =
$$
$$
= \inf\{\mu(\psi):\psi\in C_b(X), \psi\geq (\inf\{\varphi(y): y\in
Y\})_X+c_X, \psi|Y=\varphi+c_Y\}=
$$
$$=inf\{\mu(\psi-c_X)+c\cdot\mu(1_X):\psi\in C_b(X),
$$
$$
\psi-c_X\geq (inf\{\varphi(y): y\in Y\})_X, (\psi-c_X)|Y=\varphi\}=
$$
$$
=inf\{\mu(\psi-c_X):\psi\in C_b(X),
$$
$$
\psi-c_X\geq (inf\{\varphi(y): y\in Y\})_X,
(\psi-c_X)|Y=\varphi\}+c\cdot\mu(1_X)=
$$
$$
=r^X_Y(\mu)(\varphi)+c\cdot\mu(1_X)=
$$
$$
=(\mbox{by virtue of (3)})=
$$
$$=r^X_Y(\mu)(\varphi)+c\cdot r^X_Y(\mu)(1_Y),
$$
i. e. $r^X_Y(\mu)(\varphi+c_Y)=r^X_Y(\mu)(\varphi)+c\cdot
r^X_Y(\mu)(1_Y).$

Now let us show that $r^X_Y(\mu)$ is an order-preserving functional.
Let $\varphi_i\in C_b(Y),$ $i=1,2,$ and $\varphi_1\leq\varphi_2.$
Then
$$
r^X_Y(\mu)(\varphi_1)=
$$
$$
=inf\{\mu(\psi):\psi\in C_b(X),\ \psi\geq
(inf \{\varphi_1(y): y\in Y\})_X,\ \psi|Y=\varphi_1\}\leq
$$
$$
\leq inf\{\mu(\psi):\psi\in C_b(X),\ \psi\geq (inf \{\varphi_2(y):
y\in Y\})_X,\ \psi|Y=\varphi_1\}\leq
$$
$$
\leq inf\{\mu(\psi):\psi\in C_b(X),\ \psi\geq (inf \{\varphi_2(y):
y\in Y\})_X,\ \psi|Y=\varphi_2\}=
$$
$$
=r^X_Y(\mu)(\varphi_2),
$$
i. e. $r^X_Y(\mu)(\varphi_1)\leq r^X_Y(\mu)(\varphi_2).$ Thus the
functional $r^X_Y(\mu)$ is order-preserving and weakly additive on
$C_b(Y).$ Lemma 1 is proved.

The order-preserving functional $r^X_Y(\mu)$ defined as above is
said to be a \textit{restriction} of the given  order-preserving
functional $\mu\in W(\beta X)$ and the map $r^X_Y: W(\beta
X)\rightarrow W(\beta Y)$ is called the \textit{restriction
operator.}

From Lemma 1 and the equality (3) we have the following

\textbf{Proposition 1.} \textit{Let $Y$ be a $C$-embedded subspace
of Tychonoff space $X.$ Then $r^X_Y(\mu)\in O(\beta Y)$ if and only
if $\mu\in O(\beta X).$}

Let $Y\subset X,$ $\mu\in W(\beta Y)$ and $\varphi\in C_b(X).$ Put
$$
e^Y_X(\mu)(\varphi)=\mu(\varphi|Y). \eqno (4)
$$

The following statement is obvious.

\textbf{Proposition 2.} \textit{For every $\mu\in W(\beta Y)$ we
have $e^Y_X(\mu)\in W(\beta X)$ and $e^Y_X(\mu)(1_X)=\mu(1_Y).$
Hence, $e_X^Y(\mu)\in O(\beta Y)$ if and only if $\mu\in O(\beta
X).$}

The order-preserving functional $e^Y_X(\mu)$ is said to be the
\textit{extension} of the given order-preserving functional $\mu,$
and the map $e_X^Y: W(\beta X)\rightarrow W(\beta Y)$ is called the
\textit{extension operator.}

\textbf{Lemma 2.} \textit{Let $Y$ be a $C$-embedded subspace of a
Tychonoff space $X.$ Then $r^X_Y\circ e_X^Y=id_{W\circ \beta (Y)}.$}

\textsl{Proof.}  If $\mu\in W(\beta Y)$ then by virtue of
Proposition 2 $e^Y_X(\mu)\in W(\beta X)$. From Proposition 1 it
follows that $r^X_Y( e_X^Y(\mu))\in W(\beta Y).$ According to the
construction of the restriction operator the restriction $r^X_Y(
e_X^Y(\mu))$ of the order-preserving functional $e^Y_X(\mu)$ is
defined on $C=\{\psi|Y: \psi\in C_b(X)\}\equiv C_b(Y).$ But
according to (2) and (4) we have
$r^X_Y(e_X^Y(\mu))(\varphi)=\mu(\varphi)$ for each $\varphi\in
C_b(Y).$ Lemma 2 is proved.

\textbf{Lemma 3.} \textit{Let $Y$ be a $C$-embedded subspace of a
Tychonoff space $X.$ An order-preserving functional $\mu\in W(\beta
Y)$ is $\tau$-smooth if and only if $e^Y_X(\mu)\in W(\beta X)$ is a
$\tau$-smooth order-preserving functional.}

\textsl{Proof.} Let $\mu\in W_\tau(Y)$ be an arbitrary
order-preserving functional and let $\{\varphi_\alpha\}\subset
C_b(X)$ be a net such that $\varphi_\alpha\downarrow 0_X.$ Then
$\varphi_\alpha|Y\downarrow 0_Y.$ Hence,
$e^Y_X(\mu)(\varphi_\alpha)=\mu(\varphi_\alpha|Y)\rightarrow 0.$ So,
$e^Y_X(\mu)\in W_\tau(X).$

Let us establish the converse statement. Let $\mu\in W(\beta Y)$ be
an arbitrary order-preserving functional such that $e^Y_X(\mu)\in
W_\tau(X).$ Then for each net $\{\psi_\alpha\}\subset C_b(X)$
monotone decreasing to zero on $X$ we have
$e^Y_X(\mu)(\psi_\alpha)\rightarrow 0.$ Hence from (4) it follows
that for each net $\{\psi_\alpha\}\subset C_b(X)$ satisfying
$\psi_\alpha|Y\downarrow 0_Y$ one has $\mu(\psi_\alpha|Y)\rightarrow
0$.

Now take an arbitrary net $\{\varphi_\alpha\}\subset C_b(Y)$ such
that $\varphi_\alpha\downarrow 0_Y.$ Since $Y$ is $C$-embedded in
$X$ there exists a net $\{\psi_\alpha\}\subset C_b(X)$ such that
$\psi_\alpha|Y=\varphi_\alpha$ for all $\alpha.$ Hence,
$\psi_\alpha|Y\downarrow 0_Y$ and therefore
$\mu(\varphi_\alpha)=\mu(\psi_\alpha|Y)\rightarrow 0.$ Thus $\mu\in
W_\tau(Y).$ Lemma 3 is proved.

Note that for a compact $X$ each order-preserving weakly additive
functional $\mu:C(X)\rightarrow\mathbb{R}$ has a (continuous)
order-preserving weakly additive extension
$\mu':B(X)\rightarrow\mathbb{R}$ with $\mu'(1_X)=\mu(1_X)$ [5]. Here
$B(X)$ is the space of all bounded functions equipped with the
uniform convergence topology. As we have noted above for each
Tychonoff space $X$ the normed spaces $C_b(X)$ and $C(\beta X)$ are
isometrically isomorphic. Therefore any $\tau$-smooth
order-preserving functional $\mu:C_b(X)\cong C(\beta
X)\rightarrow\mathbb{R}$ may be also extended to $B(\beta X)$ as
well. We shall use the same notation for an order-preserving
functional from $W(\beta X)$ and for its extension on $B(\beta X)$.

Let $Y$ be a subspace of a Tychonoff space $X.$ Consider the
following
set\\
$O^*_Y(X)=\{\mu\in O_\tau(X): \mu(\chi_K)=0$

$\qquad \qquad \qquad \qquad \mbox{ for every compact } K\subset X
\mbox{ such that } K \cap
Y=\emptyset\}, $ \\
where $\chi_K$ is the characteristic function of the set $K.$

The equalities (2) and (4) imply the following

\textbf{Proposition 3.} \textit{Let $Y$ be a $C$-embedded subspace
of a Tychonoff space $X.$ Then $e_X^Y\circ
r^X_Y|O^*_Y(X)=id_{O^*_Y(X)}.$}

Lemmas 2, 3 and Proposition 3 yield that for a $C$-embedded subspace
$Y$ of a Tychonoff space $X$ the following equalities hold
$$
e_X^Y(O_\tau(Y))=O_Y^*(X),
$$
$$
r^X_Y(O_Y^*(X))=O_\tau(Y).
$$

These equalities imply the following

\textbf{Proposition 4.} \textit{For any Tychonoff space $X$ the maps
$$
e_{\beta X}^X: O_\tau(X)\rightarrow O^*_X(\beta X)
$$
and
$$
r^{\beta X}_X: O^*_X(\beta X)\rightarrow O_\tau(X)
$$
are mutually inverse homeomorphisms.}

The next statement is the key result.

\textbf{Theorem 1.} \textit{For an arbitrary Tychonoff space $X$ and
for every its compactification $b X$ the spaces $O^*_X(\beta X)$ and
$O^*_X(b X)$ are homeomorphic.}

\textsl{Proof.} At first recall that a continuous map $f:
b_1X\rightarrow b_2X$ between compactifications $b_1X$ and $b_2X$ of
the given Tychonoff space $X$ is called \textit{natural}, if
$f(x)=x$ for all $x\in X$ [5, P. 47].

Let $X$ be a Tychonoff space, and suppose that $b X$ is its
arbitrary compact extension. Let $f:\beta X\rightarrow bX$ be a
natural map. Assume that $\mu\in O^*_X(\beta X)$ and
$O(f)(\mu)=\nu.$ Consider an arbitrary compact set $F\subset
bX\setminus X.$ By virtue of Theorem 3.5.7 [6, P. 220] the inclusion
$f^{-1}(F)\subset \beta X\setminus X$ holds. Put $K=f^{-1}(F).$ Then
$f(K)=F$ and $\chi_K=\chi_F\circ f.$ We have
$$
\nu(\chi_F)=O(f)(\mu)(\chi_F)=\mu(\chi_F\circ f)=\mu(\chi_K)=0.
$$
So, $O(f)(O^*_X(\beta X))\subset O^*_X(b X).$ In other words the
following restriction map is correctly defined
$$
O(f)|O^*_X(\beta X):\ O^*_X(\beta X)\rightarrow O^*_X(b X). \eqno
(5)
$$
The map $O(f)|O^*_X(\beta X)$ is continuous as the restriction of
the continuous map $O(f):O(\beta X)\rightarrow O(b X).$

Let $\mu\in O(\beta X)\setminus O^*_X(\beta X).$ Then there exists a
compact set $K\subset \beta X\setminus X$ such that $\mu(\chi_K)\neq
0.$ Applying theorem 3.5.7 [6, P. 220] we obtain $f(K)\subset b
X\setminus X$ and
$$
O(f)(\mu)(\chi_{f(K)})=\mu(\chi_{f(K)}\circ f)=\mu(\chi_K)\neq 0.
$$
Hence, $O(f)(\mu)\in O(b X)\setminus O^*_X(b X).$ So
$O(f)(O^*_X(\beta X))= O^*_X(b X)$ since the map $O(f)$ is
surjective, i. e. the map  (5) is  surjective.

Note that for every compact extension $bX$ of a Tychonoff space $X$
the inclusion
$$
O_\beta(X)\subset O(bX),
$$
is true. Thus, for every Tychonoff space $X$ and its compact
extension $bX$ one has
$$
O_\beta(X)\subset O(\beta X)\cap O(bX). \eqno (6)
$$
From this it follows
$$
O_\beta(X)\subset O^*_X(\beta X)\cap O^*_X(bX).
$$
Lemma 4 [2] and the inclusion (6) imply
$$
O(f)(\mu)=\mu \eqno (7)
$$
for each order-preserving functional $\mu\in O_\beta(X).$

Now we need the density lemma for order-preserving functionals.
Recall that the \textit{density} of a topological space $X$ is the
least cardinal number of the form $|A|$ where $A$ runs over
everywhere dense subsets of the space $X,$ and $|A|$ denotes the
cardinality of the set $A.$ The density of a topological space $X$
is denoted by $d(X).$

For a Tychonoff space $X$ put
$$
O_\omega(X)=\{\mu\in O(\beta X): \mbox {supp} \mu\subset X \mbox {
and } \mbox {supp} \mu \mbox { is finite set}\}.
$$

The following statement may be considered as a version of the
density lemma 1.4 from [8] for order-preserving functionals.

\textbf{Lemma 4.} \textit{For an infinite Tychonoff space $X$ and
for its subspace $Y$ the set $O_\omega(Y)$ is everywhere dense in
$O(\beta X)$ if and only if $Y$ is everywhere dense in $X$.}

\textsl{Proof.} If $Y$ is not everywhere dense in $X$ then there
exists a nonempty open set $U\subset X$ such that $U\cap
Y=\emptyset.$ Take $x\in U.$ Consider a basic neighborhood $\langle
\delta_x; \varphi; \varphi(x)\rangle,$ where $\delta_x$ is the Dirac
measure, defined as $\delta_x(\psi)=\psi(x),$ $\psi\in C_b(X),$ and
$\varphi\in C_b(X)$ is a function such that $\varphi(x)>0$ and
$\varphi(y)=0$ for all $y\in X\setminus U.$ Then it is clear that
$\langle \delta_x; \varphi; \varphi(x)\rangle \cap
O_\omega(Y)=\emptyset.$

Let now $Y$ be an everywhere dense in $X$. Then we have
$$
O_\omega(Y)\subset O_\beta(Y)\subset \mbox{(since $O_\beta$ is
monomorphic [3])} \subset O_\beta(X),
$$
and hence, $O_\omega(Y)\subset O_\omega(X).$  Let $\mu\in
O_\omega(X)$ be an arbitrary order-preserving functional, and let
$\langle\mu;\varphi_1,...,\varphi_k;\varepsilon \rangle$ be a
neighborhood of $\mu.$ Suppose that
$\mbox{supp}\mu=\{x_1,...,x_n\}.$ One can choose a set
$\{y_1,...,y_s\}\subset Y$ and an order-preserving functional
$\nu\in O_\omega(Y)$ such that the following conditions hold:

(i) $\mbox{supp}\nu=\{y_1,...,y_n\}$;

(ii) $|\nu(\varphi_i)-\mu(\varphi_i)|<\varepsilon,$ $i=1,...,k.$\\
This implies $\nu\in \langle\mu;\varphi_1,...,\varphi_k;\varepsilon
\rangle,$ i. e. the set $O_\omega(Y)$ is everywhere dense in
$O_\omega(X).$

From the above in particular it follows that the set $O_\omega(X)$
is everywhere dense in $O_\omega(\beta X)$. On the other hand
according to proposition 3 [2] $O_\omega(\beta X)$ is everywhere
dense in $O(\beta X)$. Therefore $O_\omega(X)$ is everywhere dense
in $O(\beta X)$ and thus, $O_\omega(Y)$ is everywhere dense in
$O(\beta X).$ Lemma 4 is proved.

According to Lemma 4 the set $O_\beta(X)$ is everywhere dense in the
spaces $O(\beta X)$ and $O(bX).$ Hence, $O_\beta(X)$ is everywhere
dense in the sets $O^*_X(\beta X)$ and $O^*_X(b X)$.

Now let us show that the map (5) is one-to-one. For this purpose
take an arbitrary order-preserving functional $\nu\in O^*_X(bX).$
Suppose that there exist order-preserving functionals $\mu_1,
\mu_2\in O^*_X(\beta X)$ such that $\mu_1\neq \mu_2$ and
$O(f)(\mu_1)=O(f)(\mu_2)=\nu.$ Let $\{\mu^i_\alpha\}\subset
O_\beta(X),$ $i=1,\ 2,$ be two nets converging to the functionals
$\mu_1$ and $\mu_2,$ respectively. Since the map $O(f)|O^*_X(\beta
X):\ O^*_X(\beta X)\rightarrow O^*_X(b X)$ is continuous, the nets
$\{O(f)(\mu^i_\alpha)\},$ $i=1,\ 2,$ converge to $\nu.$ On the other
hand according to (7) one has $O(f)(\mu^i_\alpha)=\mu^i_\alpha$ for
$i=1,\ 2$ and for all $\alpha.$ Hence,
$$
\mu_1=\lim_\alpha\mu^1_\alpha=\nu=\lim_\alpha\mu^2_\alpha=\mu_2.
$$
We obtain a contradiction which shows that our assumption is false.

From the above, in particular, it follows that the map
$$
(O(f)|O^*_X(\beta X))^{-1}: O^*_X(b X)\rightarrow O^*_X(\beta X),
$$
inverse to (5), is also continuous. Thus, the map (5) is a
homeomorphism of the spaces $O^*_X(\beta X)$ and $O^*_X(b X).$
Theorem 1 is proved.

Since each Tychonoff space $X$ has a compact extension $b X$ such
that $w(X)=w(b X)$,  Proposition 4 and Theorem 1 imply the following

\textbf{Corollary 1.} \textit{The functor $O_\tau$ preserves the
weight of every infinite Tychonoff space $X$, i. e.
$w(O_\tau(X))=w(X)$.}

According to (1) we have

\textbf{Corollary 2.} \textit{The functor $O_R$ preserves the weight
of every infinite Tychonoff space $X$, i. e. $w(O_R(X))=w(X)$.}

Thus, $$w(O_\beta(X))=w(O_R(X))=w(O_\tau(X))=w(X)$$ for every
infinite Tychonoff space $X.$

Note that from Lemma 4 one can also obtain a strengthened version of
theorems 1.8 and 2.6 from [1].

Let $X$ be an infinite Tychonoff space. By virtue of the inclusions
(1) and $O_\omega(X)\subset O_\beta(X)$ it follows that $O_\beta(X)$
is everywhere dense in the spaces $O_R(X)$, $O_\tau(X)$ and $O(\beta
X)$. On the other hand according to results of [3] one has
$d(O_\beta(X))\leq d(X)$ for every infinite Tychonoff space $X.$
Thus, a strengthening of theorems 1.8 and 2.6 from [1] may be to
stated as follows

\textbf{Corollary 3.} \textit{The density of an infinite Tychonoff
space does not increase under the functors:}

\textit{$O(\beta\cdot)$ of all order-preserving functionals,}

\textit{$O_\tau$ of $\tau$-smooth order-preserving functionals, }

\textit{$O_R$ of Radon order-preserving functionals and}

\textit{$O_\beta$ of order-preserving functionals with compact
supports.}

Moreover, for every infinite Tychonoff space $X$ we have
$$
d(O(\beta X)\leq d(O_\tau(X))\leq d(O_R(X))\leq d(O_\beta (X))\leq
d(X).
$$

Recall the following notion.

\textbf{Definition  3[3].} The \textit{weak density} $wd(X)$ of a
topological space $X$ is the least cardinal number $\tau$ such that
$X$ has a $\pi$-base which is the union of $\tau$ centered families
of open sets in $X$.

We need the following properties of the weak density [3]:

\begin{tabular}{ll}
(A)&  If $Y$ is everywhere dense in  $X$ then $wd(Y)=wd(X);$\\
(B)& If $X$ is compact then $wd(X)=d(X).$
\end{tabular}

The property (A) of the weak density, Lemma 4 and the inclusions (1)
imply:

\textbf{Proposition 5.} \textit{For each Tychonoff space $X$ one
has}
$$
wd(O_\omega(X))=wd(O_\beta(X))=wd(O_R(X))=wd(O_\tau(X))=wd(O(\beta
X)). \eqno (8)
$$

Moreover according to (A) and (B) one has $wd(X)=wd(\beta X)=d(\beta
X)$ and $wd(O(\beta X))=d(O(\beta X)).$ Therefore Lemma 4 and
Proposition 5 imply

\textbf{Corollary 4.} \textit{The density of an infinite Tychonoff
space does not increase under the functors $O(\beta\cdot),$
$O_\tau$, $O_R$ and $O_\beta$.}

In connection with Proposition 5 the following question arises

\textbf{Question.} Are equalities similar to (8) valid for the
density of a Tychonoff space $X$?

The next result gives a positive answer for a particular case.

\textbf{Proposition 6.} \textit{The space $O_R(X)$ is separable if
and only if the space $O_\beta(X)$ is separable.}

\textsl{Proof.} Since $O_\beta(X)$ is everywhere dense in $O_R(X)$
one has the inequality $d(O_R(X))\leq d(O_\beta(X)).$ Let us show
that the opposite inequality is also true.

Let $\{\mu_n\}\subset O_R(X)$ be a countable everywhere dense subset
of order-preserving functionals. For every order-preserving
functional $\mu_n$ and each positive integer $m$ there exists a
compact set $K_{n, m}\subset X$ such that
$$
\mu_n(\varphi)<\frac{1}{m}\eqno (9)
$$
where $\varphi\in C_b(X)$ is an arbitrary function satisfying the
following inequalities
$$
0\leq \varphi\leq \chi_{(X\setminus K_{n,m})}.\eqno (10)
$$
Define an order-preserving functional  $\mu_{n,m}$ on $C_b(X)$ by
the formula
$$
\mu_{n,m}=r^X_{K_{n,m}}(\mu_n).\eqno (11)
$$
Then $\mu_{n,m}\in O(K_{n,m})\subset O_\beta(X).$ We have
$$|\mu_n(\varphi)-\mu_{n,m}(\varphi)| =\mbox{(according to (11))} =
|\mu_n(\varphi) - r^X_{K_{n,m}}(\mu_n)(\varphi)| =$$
$$=|\mu_n(\varphi) - r^X_{K_{n,m}}(\mu_n)(\varphi|K_{n,m})| = \mbox{(according to (10))}=$$
$$=|\mu_n(\varphi) |<\mbox{(according to (9))}< \frac{1}{m}$$
for all $n, m.$ From this it follows that the sequence
$\{\mu_{n,m}\}_{m=1}^\infty$ pointwise converges to $\mu_n$. Hence,
$M\equiv\{\mu_{n,m}:m,n=1,2,...\}$ is everywhere dense in $O_R(X).$
On the other hand, $M\subset O_\beta(X)\subset O_R(X).$ This means
that $M$ is everywhere dense in $O_\beta(X),$ i. e.
$d(O_\beta(X))\leq d(O_R(X)).$ Proposition 6 is proved.

\vspace{1cm}

\textbf{Acknowledgments.} \emph{The authors would like to
acknowledge the hospitality of the $\,$ "Institut f\"{u}r Angewandte
Mathematik",$\,$ Universit\"{a}t Bonn (Germany). This work is
supported in part by the DFG 436 USB 113/10/0-1 project (Germany)
and the Fundamental Research  Foundation of the Uzbekistan Academy
of Sciences.}

\end{document}